\documentclass{amsart}
\usepackage[T1]{fontenc}
\usepackage[utf8]{inputenc}
\usepackage[english]{babel}
\usepackage{amsthm}
\usepackage{amssymb}
\usepackage{amsmath}
\usepackage[all]{xy}
\usepackage{mathrsfs}
\usepackage{enumerate}
\usepackage{verbatim}
\usepackage{enumitem}
\usepackage{stmaryrd}
\usepackage{pgf,tikz,pgfplots}
\usetikzlibrary{arrows}
\newtheorem{teo}{Theorem}[section]
\newtheorem{thmx}{Theorem}

\newtheorem{corx}[thmx]{Corollary}

\newtheorem{lem}[teo]{Lemma}     
\newtheorem{prop}[teo]{Proposition}
\newtheorem{cor}[teo]{Corollary}

\theoremstyle{definition}
\newtheorem{defin}[teo]{Definition}
\newtheorem{ex}[teo]{Example}
\theoremstyle{claim}

\newcommand{\normal}{\trianglelefteq}

\def\normal{\trianglelefteq}

\def\fit#1{{\bf F}(#1)}
\def\fiti#1#2{{\bf F}_{#1}(#2)}

\def\zent#1{{\bf Z}(#1)}
\def\cent#1#2{{\bf C}_{#1}(#2)}
\def\oh#1#2{{\bf O}_{#1}(#2)}

\def\O#1{{\bf O}(#1)}

\def\gen#1{\langle#1\rangle}
\def\sbs{\subseteq}
\def\irr#1{{\rm Irr}(#1)}

\def\GF#1{{\rm GF}{(#1)}}
\def\field#1{{\mathbb #1}}

\def\SL#1#2{{\rm SL}_#1(#2)}

\def\syl#1#2{{\rm Syl}_{#1}(#2)}
\def\hall#1#2{{\rm Hall}_{#1}(#2)}

\def\aut#1{{\rm Aut}(#1)}
\def\out#1{{\rm Out}(#1)}

\def\gk#1{{\Gamma} (#1)}

\begin{document}
\title{On the cut-set of the Gruenberg-Kegel graph of a finite solvable group}
\author{lorenzo bonazzi}
\address{Department of Mathematics and Informatics, University of Florence}
\email{lorenzo.bonax@gmail.com}
\maketitle

\begin{abstract}
Let $\Gamma(G)$ be the Gruenberg-Kegel graph of a finite solvable group $G$. If $\sigma$ is a cut-set for $\Gamma(G)$, then $G$ has a normal $\sigma$-series of length at most $5$. As a consequence, we prove that if $G$ is a solvable group and $\Gamma(G)$ has a cut-vertex $p$, then there is a bound for the Fitting length. This bound is the best possible.
\end{abstract}
\section{Introduction}
\emph{This is a substantially revised version of \cite{bonazzi}.}\newline
Throughout this paper, all groups are assumed to be finite. If $G$ is a group, the \emph{Gruenberg-Kegel graph} $\Gamma(G)$ is defined as follows: the set of vertices consists of $\pi(G)$, the set of primes that divide $|G|$, and two vertices are joined by an (undirected) edge if and only if there is an element of $G$ of order $pq$. One of the earliest results in this topic is the well-known Gruenberg-Kegel Theorem: if $G$ is solvable and $\Gamma(G)$ is disconnected, then $G$ is a Frobenius or a $2$-Frobenius group and there are exactly two connected components, which are complete subgraphs. Refer to Theorem \ref{grukeg} a more detailed formulation.

If $\Gamma$ is a graph with vertices $V$ and $\sigma$ is a set, then $\Gamma -\sigma$ is the graph that has vertex set $V \setminus \sigma$ and whose edge relation is given as follows: two vertices of $\Gamma - \sigma $ are adjacent in $\Gamma -\sigma$ if and only if they are adjacent in $\Gamma$. We say that $\sigma $ is a \emph{psuedo cut-set} of $\Gamma$ if $\Gamma - \sigma$ is disconnected. If $\Gamma$ is connected and $\sigma \sbs V$ is a pseudo cut-set of $\Gamma$, then we call $\sigma$ a \emph{cut-set} of $\Gamma$.
If $\sigma=\{v\}$ is a cut-set for $\Gamma$, then $v$ is referred to as a cut-vertex of \emph{cut-vertex} of $\Gamma$.
The concept of a pseudo cut-set is introduced to streamline proofs, avoiding the need to distinguish cases where $\Gamma$ is disconnected. Note that if $\Gamma$ is disconnected, every $\sigma$ such that $\sigma\cap V=\emptyset$ is a pseudo cut-set for $\Gamma$. Note also that if $\sigma$ is a cut-set for $\Gamma$, then $\sigma$ is a pseudo cut-set for $\Gamma$.

We use the following notations. $\O G$ is the largest normal odd-order subgroup of $G$. If $G$ is a group, recall that the Fitting subgroup $\fit G$ is the largest normal nilpotent subgroup of $G$. If $G$ is a group we set $\fiti 1 G = \fit G$ and, for an integer $i\ge 2$, recursively  $ \fiti i G / \fiti {i-1} G = \fit {G/\fiti {i-1} G}$. If $G$ is a solvable group, there is $i$ such that $\fiti i G = G$. The smallest $i$ such that this happens is denoted $h(G)$ and it is called the \emph{Fitting length of $G$}. If $\sigma$ is a set of primes, then a $\sigma$-series of $G$ is a normal series whose factors are either $\sigma$-groups or $\sigma'$-groups. Solvable groups always have a $\sigma$-series. Denote with $\ell_\sigma(G)$, the \emph{$\sigma$-length} of $G$, is the minimal length among all the $\sigma$-series of $G$. \newline
There is a group of order 48 which plays a special role in our results, the binary octahedral group which we denote by $(2.S_4)^-$ (see Definition \ref{binary-octaedral-def}). \newline
Our main is the following refinement of the Gruenberg-Kegel theorem.

\begin{thmx}\label{theoremA}
	Let $G$ be a solvable group. Suppose that $\sigma\sbs \pi(G)$ is a pseudo cut-set for its Gruenberg-Kegel graph $\gk G$. Then, $\gk G - \sigma$ consists of two complete connected components with vertex-sets $\pi_1$, $\pi_2$ $\pi_2$, say. Subject to possibly swapping $\pi_1$ and $\pi_2$, there exists a normal series
	\[1\le G_0\le G_1\le G_2 \le G_3\le G\]
	such that $G_0$ and $G_2/G_1$ are $\sigma$-groups, $G_1/G_0$ is a nilpotent $\pi_1$-group, $G_3/G_2$ is a nilpotent $\pi_2$-group and the factor $G/G_3$ satisfies the following: either $G/G_3$ is nilpotent or  $h(G/G_3)=2$; in the latter case $2 \in \pi_2$  and $G/\O G\simeq(2.S_4)^-$.  
\end{thmx}
Theorem \ref{theoremA} has the following consequence.
\begin{corx}\label{corollaryB}	
Let $G$ be a solvable group. Suppose that $\sigma$ is a cut-set for $\Gamma(G)$, then $\ell_\sigma(G)\le 3$. Moreover, if $\sigma$ consists of a single cut-vertex $p$ of $\Gamma(G)$, then $h(G)\le 6$ and if $G/\O G\not\simeq(2.S_4)^-$, then $h(G) \le 5$.
Also, there exists a group as in the hypotheses such that $h(G) = 5$, so bounds are the best possible.
\end{corx} 
The corollary generalizes bounds for the Fitting length for solvable groups whose Gruenberg-Kegel graph is a three chain, as described in \cite[Lemma 2.3]{trianglefree}.

\section{Preliminaries}
Let $G$ be a group and $H$ a non-trivial subgroup of $G$. If $H\cap H^g =1 $ whenever $g \in G\setminus H$, then $H$ is called a \emph{Frobenius complement} in $G$. A group which contains a Frobenius complement is called a \emph{Frobenius group}. According to \cite[Theorem 2]{ctfg}, if $G$ is a Frobenius group with the Frobenius complement $H$, then there exists $N \normal G$ such that $G = N\rtimes H$. The group $N$ is called \emph{Frobenius kernel} and it is unique in $G$.
\begin{defin} If $G$ is a group, we say that $G$ is a \emph{$2$-Frobenius group} if there is a normal series $1<H<K<G$ such that $K$ is a Frobenius group with kernel $H$ and $G/H$ is a Frobenius group with kernel $K/H$. We call $H$ the \emph{lower kernel} and $K/H$ the \emph{upper kernel}. 
\end{defin}
As a consequence of Thompson's theorem on nilpotency of Frobenius kernels, we have the following.
\begin{lem} \label{frobenius-fitting}
	Let $G$ be a Frobenius group. Then $\fit G$ is the Frobenius kernel.
\end{lem}
\proof Let $L$ be the Frobenius kernel of $G$. Then $L \le \fit G$ by \cite[Theorem 10.3.1]{gorenstein}, so $\zent {\fit G} \le \cent G L\le L$ as $L$ is the Frobenius kernel of $G$. Then $\fit G \le\cent G {\zent {\fit G}} \le L$, so in fact $\fit G=L$.
\endproof

\begin{prop}\label{2frob}
	Let $G$ be a $2$-Frobenius group. Then, $\fit G$ is the lower kernel and $\fiti 2 G/\fit G$ is the upper kernel. 
	Moreover, the upper kernel is a cyclic group of odd order, $G/\fiti 2 G$ is cyclic and the lower kernel is not cyclic.
\end{prop}
\proof By definition, there is a normal series $1<H<K<G$ of $G$ such that $H$ is a Frobenius kernel of $K$ and $K/H$ is a Frobenius kernel of $G/H$. By Lemma \ref{frobenius-fitting}, we have that $\fit {G/H}=K/H$, therefore $\fit G \le H$. By Lemma \ref{frobenius-fitting} again, we have that $H=\fit K$, therefore $\fit G \le H\le \fit G$, so $H=\fit G$ and $K=\fiti 2 G$.
The remaining part of the Proposition can be found in \cite[Lemma 2.1]{primegraphsolv2015}.
\endproof 
The above Proposition is well-known; see \cite[Lemma 2]{zinoveva-mazurov-2013} for a slightly different statement.  The next proposition, is known as the "Lucido's $3$ primes lemma".
\begin{lem}\label{3primes}\cite[Proposition 1]{lucido99}
	Let $G$ be a solvable group. If $p,q,r$ are distinct primes dividing $|G|$, then $G$ contains an slement of order contained in $\{pq,pr,qr\}$.
\end{lem}
From a graph-theoretic point of view, this lemma says that that for any three distinct vertices in $\gk G$, there is at least one pair connected by an edge.\newline
The following statement of the Gruenberg-Kegel theorem is convenient for our considerations.
\begin{teo}[Gruenberg-Kegel] \label{grukeg} 
	Let $G$ be a solvable group. Suppose that $\sigma$ is a pseudo cut-set for its Gruenberg-Kegel graph $\Gamma(G)$ and let $H \in \hall {\sigma'} G$. Then $\Gamma(H)=\Gamma(G)- \sigma$ and $\Gamma(H)$ consists of two complete connected components. Moreover, one of the following holds.
	\begin{enumerate}
	    \item $H$ is a Frobenius group and the vertex-set of one component of $\Gamma(H)$ consists of the primes dividing the size of a Frobenius complement of $H$.
	    \item $H$ is a $2$-Frobenius group and the vertex-set of one component of $\Gamma(H)$ consists of the primes dividing the size of a Frobenius complement of $\fiti 2 H$.
\end{enumerate} 
\end{teo}
\proof Suppose that $\sigma$ is a pseudo cut-set of $G$. Observe that the vertex set of $\Gamma(G) -\sigma$ is equal to the vertex set of $\Gamma(H)$. In addition, two vertices of $\Gamma(H)$ are adjacent in $\Gamma(H)$ if and only if they are adjacent in $\Gamma(G)-\sigma$. This follows from the fact that every $\sigma'$-element is contained in some conjugate of $H$, since $G$ solvable and $H \in \hall {\sigma'} G$. Therefore, this implies that $\Gamma(G)-\sigma=\Gamma(H)$. By \cite[Corollary]{williams81}, $\Gamma(H)$ consists of two components and either part $1$ or part $2$ of the theorem holds, where the lower complement of $H$ is the Frobenius complement of $H$ when $H$ is a Frobenius group, and the complement of $\fiti 2 H$ when $H$ is a $2$-Frobenius group. By Lucido's three Primes Lemma, the connected components are complete subgraphs of $\Gamma(H)$.
\endproof

\section{Proof of Theorem \ref{theoremA}}
According to Theorem \ref{grukeg}, if $H$ is a group such that $\Gamma(G)$ is disconnected, then there are two connected components. The vertex-set of one of them can be identified with the primes dividing the order of the complement of either $H$ in the case that $H$ is a Frobenius group, or of $\fiti 2 H$ in the case $H$ that is a $2$-Frobenius group. Recall that if $n$ is an integer, than $\pi(n)$ denotes the set of primes that divide $n$. Therefore, we provide the following definition. 
\begin{defin} \label{defin-pi1pi2}
	Let $G$ be a solvable group and $\sigma$ be a set of primes that forms a pseudo cut-set of $\Gamma(G)$. Let $H \in \hall {\sigma'} G$. Given $H \in \hall{\sigma'}{G}$, we adopt the following definitions in light of Theorem \ref{grukeg}:
	\begin{enumerate}
		\item If $H$ is a Frobenius group, then $\pi_{2,\sigma}(G) = \pi([H:\fit H])$.
		\item If $H$ is a $2$-Frobenius group, then $\pi_{2,\sigma}(G) = \pi ([\fiti 2 H :\fit H])$.
	\end{enumerate}
	Moreover, we define, $\pi_{1,\sigma}(G)=\pi(H)\setminus\pi_{2,\sigma}(G)$. If the pseudo cut-set $\sigma$ is fixed, we write $\pi_1(G)=\pi_{1,\sigma}(G)$ and $\pi_2(G)=\pi_{2,\sigma}(G)$.
\end{defin}
Note that $\pi_1(G)$ and $\pi_2(G)$ are not standard notations. Indeed, these symbols are used in literature for labeling the connected components of $\Gamma(G)$ for a group $G$ when $\Gamma(G)$ is disconnected. Nevertheless, if $G$ is solvable, the two definitions coincide.

\begin{lem}
	Let $G$ a group and $N \normal G$. Suppose that $\sigma$ is a pseudo cut-set for $\Gamma(G)$. Then $\sigma$ is a pseudo cut-set for $\Gamma(G/N)$ and $\pi_i(G/N) \sbs \pi_i(G)$ for $i=1,2$.
\end{lem}

\begin{lem} \label{grukeg-fit-pi1-group}
	Let $G$ be a solvable group and $\sigma$ be a set of primes that is a pseudo cut-set for $\Gamma(G)$. Then $\pi_1(G)$ and $\pi_2(G)$ are the vertex-sets of the connected components of $\Gamma(G) -\sigma$ and $\fit{G/\oh \sigma G}$ is a $\pi_1(G)$-group.
\end{lem}
\proof By Theorem \ref{grukeg}, if $H \in \hall {\sigma'} G$, then $\Gamma(H)=\Gamma(G)-\sigma$ consists of two complete connected components and one of them has $\pi_2(G)$ as vertex-set. So, $\pi_1(G)$ is the vertex set of the other connected component. This is because $\pi(H)=\pi(G)\setminus\sigma=\pi_1(G)\cup \pi_2(G)$. We prove the remaining part of the lemma. Without loss of generality, we can assume that $\oh \sigma G=1$. Then $\fit G$ is a $\sigma'$-group and hence $\fit G \le H$. This means that $\fit G \le \fit H$. It follows from Theorem \ref{grukeg} that $\fit H$ is a $\pi_1(G)$-group, so $\fit G$ is a $\pi_1(G)$-group.
\endproof
The next result stands at the core of our study on Gruenberg-Kegel graphs.

\begin{prop} \label{coreresult2}
	Let $r$ be a prime, $H$ a solvable group and $V$ a faithful $\GF rH$-module. Suppose that $\sigma$ is a pseudo cut-set for $\Gamma(HV)$ and $r 	\notin \sigma$. Then, $ r \in \pi_1(HV)$ and there is $K \le H$ such that the following holds. 
	\begin{enumerate}
		\item $K$ is nilpotent and $\pi(K) \sbs \pi_2(HV)$.
		\item $\fit {H/\oh \sigma H}=K\oh \sigma H/\oh \sigma H$.
	\end{enumerate}
\end{prop}
\proof 
Let $\pi_i=\pi_i(HV)$ for $i=1,2$. Since $V$ is abelian and normal in $HV$, we have that $V \le \fit {HV}$. On the other hand, $\fit {HV} = \cent {HV} {\fit {HV}}\le \cent {HV} V =V$ since $V$ is a faithfun $H$-module. So, $V=\fit {HV}$ is an $r$-group, we have that $\oh \sigma {HV}=1$ and $r \in \pi_1$ by Lemma \ref{grukeg-fit-pi1-group}. 
By Theorem \ref{grukeg}, if $U \in \hall {\sigma'}{HV}$, we have that $\Gamma(HV)- \sigma=\Gamma(U)$ consists of two complete connected components, having vertex-sets $\pi_1$ and $\pi_2$ by Lemma \ref{grukeg-fit-pi1-group}. 
Moreover, $U$ is either a Frobenius group or a $2$-Frobenius group. Write $F=\fit U$ and note that $V=F$. Indeed, we have that $V \le F$; since $V=\fit {HV}$, the opposite inclusion also follows:
$$F =\cent U F \le \cent U V\le V.$$
Consider $\bar H=H/\oh\sigma H$. Then $\oh \sigma{\bar H}=1$ and $\fit {\bar H}$ is a $\sigma'$-group. 
Let $N$ the preimage in $H$ of $\fit {\bar H}$. Then $N/\oh \sigma H= \fit {\bar H}$. Let $U_0 \in \hall {\sigma'} H$. Since $U_0V \in \hall {\sigma'} {HV}$, possibly replacing $U$ with $U_0V$, we can assume that $U_0 \le U$. Since $V$ is an $r$ group with $r \notin \sigma$ and $U \in \hall {\sigma'} {HV}$, we have that $U=VU_0$. Write $K= U_0 \cap N$, so that $KV \normal U$ and $N=\oh \sigma H \rtimes K$. Observe that $K$ is nilpotent because $K \simeq \fit {\bar H}$. Since $V=\fit U$, it follows that $K\le\fiti 2 U$. If $U$ is a $2$-Frobenius group, we have that $V$ is the Frobenius kernel for $\fiti 2 U$ by Proposition \ref{2frob}, so $K$ is contained in a Frobenius complement of $\fiti 2 U$. If $U$ is a Frobenius group, then $V$ is the Frobenius kernel of $U$ by Lemma \ref{frobenius-fitting} and therefore $K$ is contained in a Frobenius complement of $U$. In any case, $K$ is a $\pi_2$-group by Definition \ref{defin-pi1pi2}. Note that $K\oh \sigma H/\oh \sigma H=\fit {H/\oh \sigma H}$. This concludes the proof.\endproof

\begin{lem}\label{binary-octaedral}
	Let $G$ be a solvable group. Suppose that $G/\O G$ is isomorphic to the extension of $\SL 2 3$ by a cyclic group of order 
	$2q$, where $q$ odd, and that a Sylow $2$-subgroup of $G/\O G$ is a 
	generalized quaternion group. Then 
	$G/\O G$ is isomorphic to the \emph{\texttt{SmallGroup(48,28)}}. 
\end{lem}
\proof Call $H=G/\O G$ and observe that $\O H=1$. Suppose that $H$ is the extension of $\SL 2 3$ by a cyclic group of order $2q$, with $q$ odd. Let $N \normal G$ such that $H/N$ is cyclic of order $q$. Then $N$ contains a subgroup of index $2$ that is normal in $G$ and that is isomorphic to $\SL 2 3$. Moreover, a Sylow $2$-subgroup of $N$ is a generalized quaternion group. Therefore, by direct check with the software \texttt{GAP}, up to isomorphism there is only one such a group, namely the \texttt{SmallGroup(48,28)}. Note that $\aut N=C_2\times S_4$. Let $R \in \hall {2'} H$, then $R$ acts on $N$ and $R/R_0 \lesssim \aut N$, where $R_0=\cent N R$. So, $[R:R_0]\le 3$. Take $x \in N$ of order $3$. therefore $x$ acts non trivially on $N$ and hence $\gen x \cap R_0=1$. It follows that $R=\gen x \times R_0$. Since $R_0$ is cyclic, we have that $NR \le \cent H {R_0}$ and $R_0 \le \zent H$. Thus, $R_0 \le \O H=1$.
\endproof
\begin{defin}\label{binary-octaedral-def}
	We denote by $(2.S_4)^-$ the a non-split extension of the group of order $2$
	by the group S4 such that a Sylow $2$-subgroup of this extension is isomorphic to the generalized
	quaternion group of order $16$. It has order $48$ and can also be found in GAP database as  \texttt{SmallGroup(48,28)}. Following \cite{trianglefree}, we call $(2.S_4)^-$ the \emph{binary octaedral group}.
\end{defin}

We will now prove Theorem \ref{theoremA}, which we restate here for convenience in a slightly different form.
\begin{teo} \label{normalseries}
	Let $G$ be a solvable group. Suppose that $\sigma$ is a set of primes that is a pseudo cut-set for $\Gamma(G)$. Then there is a normal series
	\[1\le G_0\le G_1\le G_2 \le G_3\le G\]
	such that $G_0$ and $G_2/G_1$ are $\sigma$-groups, $G_1/G_0$ is a nilpotent $\pi_1(G)$-group, $G_3/G_2$ is a nilpotent $\pi_2(G)$-group and $G/G_3$ is not nilpotent only if $2 \in \pi_2(G)$ and $G/\O G\simeq(2.S_4)^-$; in this case, $h(G/G_3)= 2$. 
\end{teo}  
\proof Let $\pi_i=\pi_i(G)$ for $i=1,2$. Call $G_0=\oh \sigma G$ and $\tilde G=(G/G_0)/\Phi(G/G_0)$. 
Suppose that there is a normal series $1\le \tilde G_1 \le \tilde G_2 \le \tilde G_3\le \tilde G$ such that $\tilde G_1=\fit {\tilde G}$ is a $\pi_1$-group, $\tilde G_1/\tilde G_2$ is a $\sigma$-group, $\tilde G_3/\tilde G_2$ is a nilpotent $\pi_2$-group and $\tilde G/\tilde G_3$ is not nilpotent if and only if $2 \in \pi_2$ and $\tilde G/\O{\tilde G}$ is isomorphic to $(2.S_4)^-$.
Consider $G_{i}$ the preimage of $\tilde G_i$ in $G$. Then, $ (G_1/G_0)/\Phi(G/G_0)=\fit{G/G_0}/\Phi(G/G_0)$ by \cite[III Satz 3.5]{huppertI} and it follows that $G_1/G_0=\fit {G/G_0}$, that 
is a nilpotent $\pi_1$-group. For $i\ge 2$, it is easy to see that $G_i/G_{i-1}\simeq\tilde G_i/\tilde G_{i-1}$ and $1\le G_0\le G_1 \le G_2\le G_3\le G$ satisfies the thesis of the theorem; in particular,  $G/G_3$ is not nilpotent if and only if $\tilde G/\tilde G_3$ is not nilpotent. This happens if and only if and $2 \in \pi_2$ and $\tilde G/\O{\tilde G} \simeq (2.S_4)^-$. Since $\tilde G=(G/\oh \sigma G)/\Phi(G/\oh \sigma G)$ and $\pi_2=\pi(G)\setminus (\sigma\cup \pi_1)$, we have that $\tilde G$ is the quotient of $G$ by a normal $\sigma \cup \pi_1$-group. Observe that $2 \notin \sigma \cup \pi_1$; so, we deduce that $\O{\tilde G}$ is a quotient of $\O G$ and $G/\O G\simeq (2.S_4)^-$. We proved that if $\tilde G$ possesses a series that satisfies the thesis of the Theorem, then the same is true for $G$. Hence, it is no loss to assume $\oh \sigma G=1$ and $\Phi(G)=1$.\\
Let $\pi_i = \pi_i(G)$ for $i=1,2$. By Lemma \ref{grukeg-fit-pi1-group}, we have that $\fit G$ is a $\pi_1$-group. 
By Gasch\"{u}tz's Theorem \cite[1.12]{mw}, $\fit G$ has a complement $H$ in $G$ and $\fit G$ is a faithful completely reducible $H$-module, possibly of mixed characteristic. Write $\fit G=M_1 \times \dots \times M_n$ as the product of irreducible $H$-modules, so that $M_i$ is an elementary abelian $r_i$-group, with $r_i \in \pi_1$ for every $i$.\\
Call $H_i=H/\cent H {M_i}$ and $\bar H=\prod H_i$. Note that $H \lesssim \bar H$, since $\bigcap_i \cent H {M_i}=\cent H {\fit G}=1$. 
The group $M_i$ is a faithful irreducible $H_i$-module. Note that $M_iH_i=G/\cent H{M_i}\prod_{j\neq i}M_j$, so $\Gamma(M_iH_i)$ is a subgraph of $\Gamma(G)$ for every $i$. Let $L \in \hall {\pi_2} H$; since no vertex in $\pi_2$ is adjacent in $\Gamma(G)$ to any vertex in $\pi_1$, we have that $\cent L {M_i}=1$. Therefore, for every $i$, $H_i$ contains a subgroup that is isomorphic to $L$. In particular, $\pi_2 \subseteq \pi(H_i)$. Since $r_i \in\pi_1$, $r_i\in  \pi(H_iM_i)$ and $\Gamma(H_iM_i) \le \Gamma(G)$, it follows that $\Gamma(M_iH_i)-\sigma$ is disconnected, for every $i$. This means that $\sigma$ is a pseudo cut-set for $\Gamma(M_iH_i)$. Note furthermore that $\pi_j(M_iH_i) \sbs \pi_j$, for $j=1,2$ and for all $i$.
By Proposition \ref{coreresult2}, for every $i$, there are $H_{i,2},H_{i,3} \normal H_i$ such that $H_{i,2}= \oh \sigma{H_i}$ and $\fit{H_i/H_{i,2}}=H_{i,3}/H_{i,2}$. Moreover, we have that $H_{i,3}=K_iH_{i,2}$, where $K_i \in \hall {\pi_2} {H_{i,3}}$ and $K_i$ is nilpotent. Suppose that there is $i$ such that $H_i/H_{i,3}$ is not nilpotent. Note that $K_i$ is nilpotent and acts regularly on $M_i$; so, by \cite[Theorem 10.3.1]{gorenstein}, we have that $K_i=C_i\times D_i$ where $C_i$ is a cyclic group of odd order and $D_i$ is a $2$-group that is cyclic, quaternion or generalized quaternion. 
Moreover,  $H_i/H_{i,3} \lesssim \out {K_i}=\out {C_i}\times \out {D_i}$, where $\out {C_i}$ is abelian. Hence, since $H_i/H_{i,3}$ is not nilpotent, $D_i$ is the quaternion group of order $8$ and $H_i/H_{i,3}$ acts on $D_i$ as the symmetric group $S_3$ (and hence $2$ divides $[H_i:H_{i,3}]$). In this case, $h(H_i/H_{i,3})=2$. In particular, $2 \in \pi_2$ and therefore that $\cent G {M_i}$ has odd order. Let $T \in \syl 2 H$; observe that $T \in \syl 2 G$, that $T \lesssim H_i$ and that the image of $T$ in $H_i$ acts regularly on $M_i$. It follows that $T$ is cyclic or a generalized quaternion group by \cite[Theorem 10.3.1]{gorenstein}; since $D_i$ is the quaternion group of order $8$ and $2$ divides $[H_i:H_{i,3}]$, we have that $T$ is the generalized quaternion of order $\ge 16$. Since $G/\O G$ is solvable, by \cite[Proposition 2]{trianglefree} we have that either $G/\O G\simeq T$ or $G/\O G$ is isomorphic to an extension of the group $\SL 2 3$ by a cyclic group of order $q$ or $2q$ with $q$ odd. Since $H/H_{i,3}$ acts on $D_i$ as the symmetric group $S_3$, the group $G$ is not $2$-nilpotent. If $G/\O G$ is isomorphic to an extension of the group $\SL 2 3$ by a cyclic group whose order $q$ with $q$ odd, then $|T|=8$, $T \simeq D_i$ and $H_i/G_2^i$ acts on $D_i$ as a cyclic group of order $3$, in this case $h(H/H_{i,3})=1$, against our assumptions. Therefore, by Lemma \ref{binary-octaedral},  $G/\O G\simeq(2.S_4)^-$.\\
Call $\bar H_2= \prod_i H_{i,2}$, $\bar H_3= \prod_i H_{i,3}$, so $1\le \bar H_2 \le \bar H_3 \le \bar H$ is a normal series such that $\bar H_2$ is a $\sigma$-group and $\bar H_3/\bar H_2$ is a nilpotent $\pi_2$-group.
Recall that there is a monomorphism $i\colon H \to \bar H$ that maps $h\mapsto \bar h$. Write 

\begin{align*}
	\phi \colon G &\to F\rtimes \bar H\\
	fh &\mapsto f\bar h.
\end{align*}
Observe that $\phi(f_1h_1f_2h_2) =f_1(f_2^{h_1^{-1}})\bar h_1\bar h_2= f_1f_2^{{\bar h_1}^{-1}}\bar h_1\bar h_2 = \phi(f_1h_1)\phi(f_2h_2)$, so $\phi$ is a morphism. Clearly, $\phi$ is a monomorphism. Take $G_1 = \fit G$, $G_2 = \phi^{-1}(\fit G (i(H)\cap \bar H_2))$ and $G_3 = \phi^{-1}(\fit G (i(H)\cap \bar H_3))$. Note that $G_1$ is a $\pi_1$-group, $G_2/G_1$ is a $\sigma$-group and $G_3/G_2$ nilpotent $\pi_2$-group. Moreover $G/G_3 \simeq H/H\cap \bar H_3 \simeq H\bar H_3/\bar H_3$ is not nilpotent if and only if there is one $i$ such that $H_i/H_{i,3}$ is not nilpotent and $2 \in \pi_2$. This happens if and only if $G/\O G \simeq (2.S_4)^-$ and in this case, $h(G/G_3)=2$.
\endproof

If $G$ is a group and $\sigma$ is a set of primes, a \emph{$\sigma$-series} is a normal series $1=G_0<G_1<\dots < G_n=G$ such that $G_i/G_{i-1}$ is a $\sigma$-group or a $\sigma'$-group for every $i=1\dots n$. The \emph{$\sigma$-length} $\ell_{\sigma}(G)$ of a group $G$ is the minimum possible number of factors that are $\sigma$-groups in any $\sigma$-series. The notion of $\sigma$-length is well-behaved with subgroups and quotient.
Moreover, if $h(G)=n$, then $\ell_\sigma(G) \le n-1$.
Now we restate Corollary \ref{corollaryB} . 
\begin{cor} \label{coreresult3}
	Let $G$ be a solvable group. Suppose that $\sigma$ is a cut-set for $\Gamma(G)$, then $\ell_\sigma(G)\le 3$. Moreover, if $\sigma$ consists of a cut-vertex $p$ of $\Gamma(G)$, then $h(G)\le 6$ and if $G/\O G\not\simeq(2.S_4)^-$, then $h(G) \le 5$.
	The bounds are the best possible.
\end{cor}

\proof By Theorem \ref{normalseries}, there is a series $1 \le G_0 \le G_1\le G_2 \le G_3 \le G$ such that $G_0$ and $G_2/G_1$ are $\sigma$-groups, $G_1/G_0$ is a $\pi_1(G)$-group and $G_3/G_2$ is a $\pi_2(G)$-group. Now, $h(G/G_3)\le 2$ and therefore $\ell_\sigma (G/G_3)\le 1$. 
Thus, $\ell_\sigma(G)\le 3$.\newline
Suppose now that $\sigma=\{p\}$. Clearly $G_0$ and $G_2/G_1$ are nilpotent because they are $p$-groups. The factor $G/G_3$ is not nilpotent if and only if $G/\O G\simeq(2.S_4)^-$. In this case, $h(G/G_3)=2$. So $h(G) \le 5$ except when $h(G/G_3)=2$. In this case, $G/\O G\simeq(2.S_4)^-$ and $h(G) \le 6$. \newline
For proving that the bounds obtained are the best possible, it suffices to assume that $\sigma=\{p\}$, where $p$ is a cut-vertex. 
If $G$ is the group constructed in \cite[Remark 2]{trianglefree}, then $G$ is solvable, $3$ is a cut-vertex of $\Gamma(G)$ and $\ell_{\sigma}(G)=3$, where $\sigma=\{3\}$. Moreover, $h(G) =6$. Observe that $G/\O G\simeq (2.S_4)^-$.  Suppose that $G/\O G\not\simeq(2.S_4)^-$, then Theorem \ref{cut-vertex-length-5} below provides an example of a group $G$ of odd order such that $h(G)=5$ and $\Gamma(G)$ has a cut-vertex.
\endproof

\section{Examples}\label{section-examples}
\begin{ex}
	Let $G$ be a solvable group, and let $\sigma$ be a cut-set for $G$. If $|\sigma|=1$, then there is an upper bound for $h(G)$ as indicated by Corollary \ref{coreresult2}. However, if $|\sigma|\ge 2$, there exists no such bound for $h(G)$.  In fact, let $n\ge 2$ be a large integer. It is not difficult to find a group $G_1$ of order $p^\alpha q^\beta$ with Fitting length $n$, where $p,q\ge 5$ are two primes. Consider $G_2=S_3$, so $\Gamma(G_2)$ is the union of two connected components that consist of the prime $2$ and the prime $3$. Let $G=G_1\times G_2$, then $\{p,q\}$ is a cut-set for $\Gamma(G)$ and  $h(G)= n$.
\end{ex}
In the hypotheses of Corollary \ref{coreresult3}, if $G/\O G\not \simeq (2.S_4)^-$, then the bound obtained is the best possible. In fact, in Theorem \ref{cut-vertex-length-5} we construct a group $G$, of odd order and arbitrarily large derived length, such that $\Gamma(G)$ has a cut-vertex and $h(G)=5$. \newline
Some concepts of representation theory are involved; we adopt the notation in \cite[Ch. 9]{ctfg}. Consider $r$ be a prime, $\tilde R$ be the ring of local integers for the prime $r$ (see \cite[pag. 265]{ctfg}). Let $ * :  \tilde R \to \tilde R/M$ be the projection map of $\tilde R$ onto the quotient $\tilde R/M$, where $M$ is a maximal ideal of $\tilde R$ containing $r$. Throughout this section, $\field F$ denotes the field $\tilde R/M$.

\begin{lem} \label{principalcomponent}
	Let $G$ be a group and $ x \in G$. Let $\chi$ be the ordinary character of $G$ afforded by the representation $\mathfrak{X}\colon G \to \text{\rm GL}(V)$. Then $[\chi_{\langle x \rangle}, 1_{\langle x \rangle}]=0$ if and only if $\mathfrak{X}(x)$ acts fixed-point-freely on $V$.
\end{lem}
\proof The principal character $1_{\langle x \rangle}$ appears among the irreducible constituents of $\chi_{\langle x \rangle}$ if and only if $\mathfrak{X}(x)$ has one eigenvector with eigenvalue $1$, namely there is a fixed point.\endproof

Let $\mathfrak X$ a complex representation of a group $G$ and suppose that, for every $g \in G$, $\mathfrak X(g)$ has entries in $\tilde R$. Then, following \cite[pag.  266]{ctfg}, if $\field F=\tilde R/M$, we can construct an $\field F$-representation $\mathfrak X^*$ of $G$ by setting $\mathfrak X^*(g)=\mathfrak X (g)^*$, where $\mathfrak X(g)^*$ is the matrix obtained by  applying $*$ to every entry of $\mathfrak X (g)$.
If $\field E \sbs \field F$ is a subfield and $\mathfrak Z$ is an $\field E$-representation of $G$, then $\mathfrak Z$ maps $G$ into a group of non-singular matrices over $\field E$. Thus, we can view $\mathfrak Z$ as an $\field F$-representation of $G$. As such we denote it by $\mathfrak Z^{\field F}$ (see \cite[pag. 144]{ctfg}). If $\mathfrak{X}$, $\mathfrak{Y}$ are two $G$-representations that are similar over some field, we write $\mathfrak{X} \simeq \mathfrak{Y}$.
\begin{lem} \label{coprimechar}
	Let $\mathfrak{X}$ be an irreducible $\mathbb{C}$-representation of a group $G$. Suppose $r$ is a prime such that $r\nmid |H|$ and $\field F = \tilde R/M$, where $\tilde R$ is the ring of local integers and $M$ a maximal ideal of $\tilde R$ containing $r$. Then, there is a finite field $\field E \subseteq \field F$, a  $\mathbb{C}$-representation $\mathfrak{Y}$ similar to $\mathfrak{X}$ that takes values in $\tilde R$ and an absolutely irreducible $\field E$-representation $\mathfrak{Z}$ such that $\mathfrak{Y}^*\simeq \mathfrak{Z}^{\field F}$.
\end{lem}
\proof
Let $\mathfrak{X}$ be a $\mathbb{C}$-representation of a group $H$ and $\chi$ be its complex character. By \cite[Theorem 15.8]{ctfg}, there exists a $\mathbb{C}$-representation $\mathfrak{Y}$, similar to $\mathfrak{X}$, that takes values in $\tilde R$, namely $\mathfrak{Y}(g) \in M_{\chi(1)}(\tilde R)$ for all $g\in G$. Moreover $\mathfrak{Y}^*$ is an $\field F$-representation of $G$ and $\hat \chi$ is its Brauer character. Since $r \nmid |G|$, by \cite[Theorem 15.13]{ctfg} we have $\hat \chi =\chi$ and $\hat \chi$ is irreducible. Hence $\mathfrak{Y}^*$ is irreducible. The field $\field F$ is algebraically closed over its prime field $\field F_p$ by \cite[Lemma 15.1c)]{ctfg}, so $\mathfrak{Y}^*$ is absolutely irreducible by \cite[Corollary 9.4]{ctfg}. Let $\field E \subset \field F$ a splitting field for the polynomial $x^{|G|}-1 \in \field F_p[x]$. Note that $\field E$ is a finite-degree extension of the prime field of $\field F$, therefore $\field E$ is finite. For every $g \in G$, $\chi^*(g)$ is a sum of $|G|$-roots of unity, so $\chi^*(g) \in \field E$ for every $g \in G$. By \cite[Theorem 9.14]{ctfg} there exists an absolutely irreducible $\field E$-representation $\mathfrak{Z}$ of $G$ such that $\mathfrak{Z}^{\field F}\simeq \mathfrak{Y}^*$.
\endproof
\begin{prop} \label{coprimechar1}
	Let $\mathfrak{X}$ be a $\mathbb{C}$-representation for a group $G$ and $W$ the associated $\mathbb{C}[G]$-module. Then, there is a finite splitting field $\field E$ of $G$ and a $\field E[G]$-module $V$, such that $(|G|,|V|)=1$ and $$\dim_{\mathbb{C}}\cent W x=\dim_{\mathbb{E}}\cent V x$$ for every $x \in G$.
\end{prop}
\proof Let $r$ be a prime that does not divide $|G|$ and denote $\tilde R$ the ring of local integers at the prime $r$. By Lemma \ref{coprimechar}, there is a $\mathbb{C}$-representation $\mathfrak{Y}$, similar to $\mathfrak{X}$, that take values in $\tilde R$, a finite field $\field E \subseteq \field F=\tilde R/M$ (where $M$ a maximal ideal of $\tilde R$ containing $r$) and absolutely irreducible $\field E$-representation $\mathfrak{Z}$ such that $\mathfrak{Y}^* \simeq \mathfrak{Z}^{\field  F}$. Call $W$ the $\mathbb{C}[G]$-module associated to $\mathfrak{X}$ and $V$ the $\field E[G]$-module associated to $\mathfrak{Z}$. Note that $V$ is finite since $\field E$ is finite, moreover $(|G|,|V|)=1$. Let $x \in G$, the number $m=\dim_{\mathbb{C}}\cent W x$ is the geometric multiplicity of the eigenvalue $1$ of the matrix $\mathfrak{Y}(x)$, that is equal to the algebraic multiplicity of $1$ of the matrix $\mathfrak{Y}(x)$. This is because the characteristic of $W$ is $0$ and thus the action of $\gen x$ on  $W$ is completely reducible by Maschke's Theorem. Hence, following the proof of \cite[Lemma 2.15]{ctfg}, $\mathfrak Y(x)$ is diagonalizable and, therefore, algebraic and geometric multiplicities of $\mathfrak Y (x)$ coincide. Note that $m$, as the algebraic multiplicity of $1$ in $\mathfrak Y(x)$, is equal to the algebraic multiplicity of the $1$ in $\mathfrak{Y}^*(x)$. Since $\mathfrak{Y}^*(x)$ and $\mathfrak{Z}^{\field F}(x)$ are similar, $m$ is the algebraic multiplicity of $1$ for the matrix $\mathfrak Z^{\field F}(x)$, that is the algebraic multiplicity of $1$ in $\mathfrak{Z}(x)$. 
Using the same argument as above, since the characteristic of $V$ does not divide $|G|$, $m$ is the geometric multiplicity of $1$ for $\mathfrak Z (x)$, that is equal  to $\dim_{\mathbb{E}}\cent V x$, 
\endproof

\begin{figure}[h!]
	\caption{}
	\label{figure}
	\includegraphics[width=0.5\textwidth]{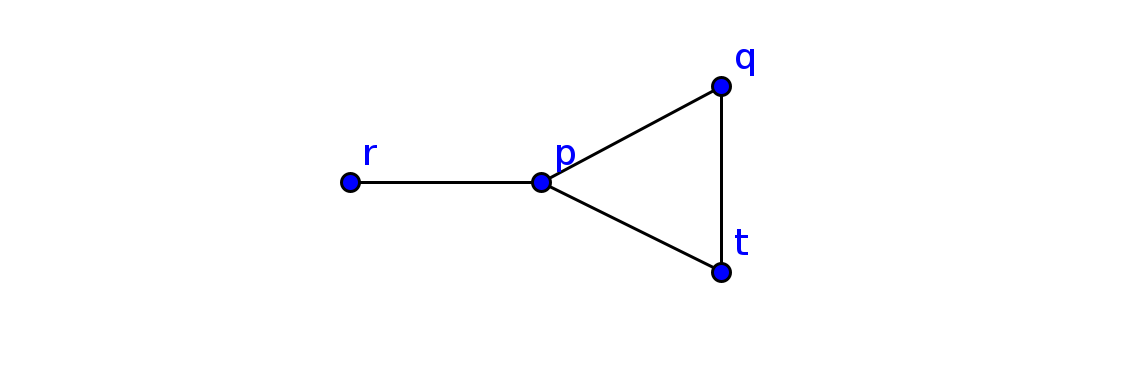}
\end{figure}

Now we prove that the second part of Corollary \ref{corollaryB} holds, specifically that the bounds we obtained in Corollary \ref{coreresult3} are the best possible.
\begin{teo} \label{cut-vertex-length-5}
	Let $n\ge5$ an integer. Then, there is a solvable group $G$ of odd order, with derived length greater than $n$, such that $h(G)=5$ and $\Gamma(G)$ is the graph in Figure \ref{figure}.
\end{teo}
\proof
Let $t,q$ be two odd primes and $T,Q$ two cyclic groups of order respectively $t^2$ and $q^2$. Let $T_0 \le T$ the subgroup of order $t$ and $Q_0\le Q$ the subgroup of order $q$. Choosing $t$ to be a prime divisor of $q^2-1$, then $t$ divides $|\aut  Q|$, and so there is a monomorphism $i\colon T/T_0 \to \aut Q$. consider $j\colon T \to T/T_0$ the projection map. Then the composition $i\cdot j$ is an action of $T$ on $Q$ whose kernel is $T_0$. Consider $L=Q\rtimes T$, $L_0=Q_0T_0$ and $\bar L=L/L_0$. Let $p$ be an odd prime. By \cite[Theorem 22.25]{berkovich}, there is a $p$-group $P$ of derived length $n$ that has a faithful irreducible character $\theta$. If $P_0$ is the base group of $P\wr \bar L$, there is an action of $L$ on $P_0$ and the kernel of such an action is $L_0$. Call $H=P_0 \rtimes L$, note that $\fit H=P_0L_0$, $\fiti 2 H=P_0F(L)$ and $\fiti 3 H=H$. Moreover, the derived length of $H$ is greater than $n$. Now,
$ P_0=\prod_{j \in \bar L} P_j $
with $P_j \simeq P$ for every $j$ and there is a character $\theta^j$ of $P_j$ that is isomorphic to $\theta$. Consider 
$\psi=\prod_j \theta^j $.
By construction, $\psi$ is a faithful irreducible character of $P_0$. Consider now two non-trivial characters  $\lambda \in \irr {Q_0}$ and $\mu \in \irr{T_0}$ such that $\lambda \mu$ is a faithful irreducible character of $L_0$. Note that $L_0P_0=P_0\times Q_0 \times T_0$, hence $\psi\lambda\mu \in \irr{P_0L_0}$. Let $\chi \in \irr{H\mid\psi\lambda\mu}$, it is easy to see that $\chi$ is faithful. Indeed $\ker\chi \cap P_0=1$ since $[\chi_{P_0}, \psi]\neq 0$ and $\psi$ is faithful. Moreover, if $q$ divides $|\ker \chi|$, we have that $Q_0 \le \ker\chi$, but $[\chi_{Q_0},\lambda]\neq 0$ and this is impossible. Replacing $q$ by $t$, we have that $t$ does not divide $|\ker \chi|$. So, we have that $\ker(\chi)=1$ and $\chi$ is faithful.
The same argument also implies that $[\chi_{Q_0},1_{Q_0}]=[\chi_{T_0},1_{T_0}]=0$. Therefore, by Lemma \ref{principalcomponent}, if $\mathfrak{X}$ is a representation fo $H$ that affords $\chi$ and $x$ is an element of order $t$ or $q$, then $x$ is contained in of either $T_0$ or $Q_0$ and $\mathfrak X (x)$ acts fixed point freely on $W$, the $\mathbb{C}[H]$-module associated to $\chi$. Let $r$ be an odd prime such that $r \nmid |H|$. By Proposition \ref{coprimechar1} there is a finite field $\field E$ of characteristic $r$ and a finite $\field E[H]$-module $V$ such that, for every $x \in H$ \[\dim_{\mathbb{C}}\cent W x=\dim_{\mathbb{E}}\cent V x.\]
Note that, since $W$ is faithful, $V$ is faithful. Moreover, an element $x \in H$ acts fixed-point-freely on $V$ whenever it does on $W$. So, every element in $H$ of order $t$ or $q$ acts fixed-point-freely on $V$ and $HV$ does not contain any element of order $tr$ or $qr$.
\newline
On the other hand, the subgroup $P_0$ has an elementary abelian subgroup of order $p^2$. Hence the action of $P_0$ on $V$ is not regular by \cite[Theorem 10.3.1]{gorenstein}. So, there is one element of order $rp$. In addition, in $HV$ there are elements of order $tp$, $qp$ and $tq$ since $\fit H=P_0L_0$. This means that $p$ is a cut-vertex for the connected graph $\Gamma(HV)$. Note that $\fit {HV}=V$ and $HV$ has Fitting length $4$. Now consider $C_p$ a group of order $p$ and let $G=C_p\wr (HV)$. Note that $G$ has odd order, $h(G)=5$ and that $G$ has derived length greater than $n$. \endproof

As per the reviewer's suggestion, groups $G$ for which $\Gamma(G)$ is isomorphic to the graph in Figure \ref{figure} are discussed in \cite{grukeg-paw-solvable}. Notably, the group constructed in the aforementioned Theorem corresponds to case 1$a$) as delineated in Theorem 2 of that publication (which was not accessible at the time of composition).

\end{document}